\documentclass [a4paper,12pt]{amsart}
\usepackage{graphicx}
\usepackage[ansinew]{inputenc}    
\usepackage[all]{xy}

\usepackage{amssymb,amscd}

\newtheorem{thm}{Theorem}[section]

\def\dim{\operatorname{dim}}

\def\fim {{\hfill $\blacksquare$}} 
\usepackage{graphicx}
\usepackage[ansinew]{inputenc}    
\usepackage[all]{xy}
\usepackage{amssymb,amscd}
\usepackage{color}

\newtheorem{cor}[thm]{Corollary}
\newtheorem{teo}[thm]{Theorem}
\newtheorem{lem}[thm]{Lemma}

\theoremstyle{definition}
\newtheorem{ex}[thm]{Example}
\newtheorem{defi}[thm]{Definition}

\def\dim{\operatorname{dim}}

\begin{document}
\title{$\mu$-constant deformations of functions on ICIS}

\author{R. S. Carvalho, B. Oréfice-Okamoto, J. N. Tomazella}

\address{Departamento de Matem\'atica, Universidade Federal de S\~ao Carlos, Caixa Postal 676,
13560-905, S\~ao Carlos, SP, BRAZIL}

\email{rafaelasoares@dm.ufscar.br}

\address{Departamento de Matem\'atica, Universidade Federal de S\~ao Carlos, Caixa Postal 676,
13560-905, S\~ao Carlos, SP, BRAZIL}

\email{bruna@dm.ufscar.br}

\address{Departamento de Matem\'atica, Universidade Federal de S\~ao Carlos, Caixa Postal 676,
13560-905, S\~ao Carlos, SP, BRAZIL}

\email{tomazella@dm.ufscar.br}

\thanks{The first author was partially supported by CAPES.
The second author was partially supported by FAPESP Grant
2016/25730-0. The third author was partially supported by CNPq Grant
309086/2017-5 and FAPESP Grant 2016/04740-7.}

\begin{abstract}
We study deformations of holomorphic function germs $f:(X,0)
\rightarrow \mathbb{C}$ where $(X,0)$ is an ICIS. We present
conditions for these deformations to have constant Milnor number,
Euler obstruction and Bruce-Roberts number.

\end{abstract}

\maketitle
\section{Introduction}

Let $f:(\mathbb{C}^n,0)\rightarrow (\mathbb{C},0)$ be a holomorphic
function germ with isolated singularity and let
$F:(\mathbb{C}^n\times\mathbb{C},0)\rightarrow (\mathbb{C},0)$ be a
deformation of $f$.

We denote, for each $t\in \mathbb{C}$,
$f_t:(\mathbb{C}^n,0)\rightarrow (\mathbb{C},0)$ the germ defined by
$f_t(x)=F(x,t)$.

Then, we have a family of function germs, $f_t$. Many authors
studied the properties of such a family and a very important result
is to know when the family has constant topological type. In this
direction, we have the Milnor number (\cite{Milnor}), which is a
well known number related to a function germ. We know that a family
$f_t$ has constant topological type if and only if, it has constant
Milnor number (\cite{Le Ramanujam}, \cite{Timourian}).

However, the problem of determining if a family has constant Milnor
number is not easy. Greuel \cite[p.161]{Greuel} presents algebraic
methods to verify if a family has constant Milnor number in the
following theorem.

\begin{teo} \cite[Theorem 1.1]{Greuel} \label{teorema 1.1 greuel} Let $f:(\mathbb{C}^n,0)\rightarrow (\mathbb{C},0)$ be a holomorphic
function germ with isolated singularity at the origin. For any
deformation $F:(\mathbb{C}^n\times \mathbb{C})\rightarrow
(\mathbb{C},0)$ of $f$ the following statements are equivalent:

(1) $F$ is $\mu$-constant.

(2) For every holomorphic curve $\gamma: (\mathbb{C},0)\rightarrow
(\mathbb{C}^n\times \mathbb{C},0)$

\begin{center}$\nu\left((\partial F/\partial t)\circ\gamma\right)>inf\{\nu\left((\partial F/\partial x_i)\circ\gamma\right) | \
    i=1,\ldots,n\}$,\end{center}

(where $\nu$ denotes the usual valuation of a complex curve).

(3) Same statement as in (2) with $``>"$ replaced by $``\geq"$.

\vspace{0.3cm}

(4) $\partial F/\partial t \in \overline{J}$, (where $\overline{J}$
denotes the integral closure of the Jacobian ideal $J=\left\langle
\partial F/\partial x_1,\ldots,\partial F/\partial x_n
\right\rangle$ in $\mathcal{O}_{n+1}$).

\vspace{0.3cm}

(5) $\partial F/\partial t \in \sqrt{J}$, (where $\sqrt{J}$ denotes
the radical of $J$).

\vspace{0.3cm}

(6) $v(J)=\{(x,t) \in \mathbb{C}^n\times\mathbb{C} \ | \ (\partial
F/\partial x_i)(x,t)=0, \ i=1,\ldots,n\}=\{0\}\times \mathbb{C}$
near $(0,0)$.

\end{teo}

Given an analytic variety germ $(X,0)$ and a function germ
$f:(\mathbb{C}^n,0)\rightarrow(\mathbb{C},0)$ we have the
Bruce-Roberts number, $\mu_{BR}(X,f)$ (\cite{Bruce e Roberts}). The
Bruce-Roberts number generalizes the Milnor number in the sense that
the germ is $\mathcal{R}_X$-finitely determined if and only if the
Bruce-Roberts number is finite, where $\mathcal{R}_X$ is the group
of the diffeomorphisms which preserves $(X,0)$.

In \cite{Tomazella}, Ahmed, Ruas and Tomazella study the Theorem
\ref{teorema 1.1 greuel} for the Bruce-Roberts number.

In \cite{Hamm 2}, Hamm introduces the Milnor number of an isolated
complete intersection singularity (ICIS). The constancy of the
topological type implies the constancy of the Milnor number of a
family of ICIS of dimension $d\neq 2$ (see \cite{Parameswaran}).

If $(X,0)\subset(\mathbb{C}^n,0)$ is an ICIS and $f:(X,0)\rightarrow
(\mathbb{C},0)$ is a holomorphic function germ with isolated
singularity, we can, also, look to the Milnor number of
$f:(X,0)\rightarrow (\mathbb{C},0)$, $\mu(f|_X)$. By the Lê-Greuel
formula (see \cite[p.77]{Looijenga})

\begin{center} $\mu(f|_X)=\mu(X,0)+\mu(X\cap f^{-1}(0),0)$,\end{center} where $\mu(X,0)$ and $\mu(X\cap f^{-1}(0),0)$ denote the Milnor
number of the ICIS defined by Hamm. Therefore, to study the
constancy of $\mu(X\cap f_t^{-1}(0))$ is equivalent to study the
constancy of $\mu(f_t|_X)$.

In this work, we study the constancy of the Milnor number of such a
family, $f_t:(X,0)\rightarrow(\mathbb{C},0)$. That is, we analyze
the assertions of Theorem \ref{teorema 1.1 greuel} for this singular
case.

We apply this study to produce examples of families of function in
an ICIS which have constant Milnor number. Moreover, we present
conditions of the constancy of the local Euler obstruction and of
the Bruce-Roberts number. Furthermore, we present a sufficient
condition for family $f_t$ to be $C^0-\mathcal{R}_X$-trivial.

Finally, we analyze the constancy of Milnor number of a family
$f_t:(X_t,0)\rightarrow(\mathbb{C},0)$ where $(X_t,0)$ is a
deformation of ICIS $(X,0)$. For this study we use the strict
integral closure of a module and we analyze the assertions of
Theorem \ref{teorema 1.1 greuel} which make sense in this new
context.

\section{Preliminary concepts}

Let $(X,0)\subset(\mathbb{C}^n,0)$ be the ICIS defined by the zero
set of a holomorphic map germ $\phi:(\mathbb{C}^n,0)\rightarrow
    (\mathbb{C}^p,0)$, with $n>p$.

Let
$\Phi:(\mathbb{C}\times\mathbb{C}^n,0)\rightarrow(\mathbb{C}^p,0)$
be a holomorphic deformation of $\phi$, defined by
$\Phi(s,z)=\phi_s(z)$, where $\phi_0=\phi$ and $X_s:=\phi_s^{-1}(0)$
is smooth for $s\neq 0$ small enough. We denote
$\mathcal{X}=\Phi^{-1}(0)$. Let $f:(X,0)\rightarrow
    (\mathbb{C},0)$ be a
holomorphic function germ and let
$F:(\mathcal{X},0)\rightarrow(\mathbb{C},0)$ such that for all
$s\neq 0$ small enough, the germ

\begin{center} $\begin{array}{cccc}
  f_s: & (X_s,0) & \rightarrow & (\mathbb{C},0) \\
   & z & \mapsto & F(s,z)
\end{array}$ \end{center} is a Morse function. Inspired by \cite{Milnor}, we define the
    Milnor number of $f$ by

    \begin{center} $\mu(f|_X)=\sharp S(f_s)$, \end{center} where $f_s:(X_s,0)\rightarrow
    (\mathbb{C},0)$ is a morsification of $f:(X,0)\rightarrow
    (\mathbb{C},0)$, with $(X_s,0)$ a smoothing of $(X,0)$  and $S(f_s)$ denotes the set of singular points
    of $f_s$.

In order to calculate this number, let $D=\{(s,z)\in\mathcal{X} \ |
\ z$ is a singular point of $f_s\}$ and $\pi:D\rightarrow \mathbb
{C}$ the restriction of the projection on the first coordinate. We
have that

\begin{center} $\mu(f|_X)=\sharp S(f_s)=$degree$(\pi)$. \end{center}

Therefore, from \cite{Mumford},

\begin{center} $\mu(f|_X)=e\left(\langle
s\rangle,\frac{\mathcal{O}_{\mathcal{X},0}}{J(f_s,\phi_s)}\right)$,
\end{center} where $\mathcal{O}_{\mathcal{X},0}$ is the local ring of
$(\mathcal{X},0)$, $J(f_s,\phi_s)$ is the ideal generated by the
maximal order minors of the Jacobian matrix of $(f_s,\phi_s)$
(partial derivatives with respect to $z$, only) and $e(I,R)$ denotes
the Hilbert-Samuel multiplicity of the ideal I with respect to the
ring R.

Since $\frac{\mathcal{O}_{\mathcal{X},0}}{J(f_s,\phi_s)}$ is a
determinantal ring, it is Cohen-Macaulay, therefore, from
\cite[p.138]{Matsumura}

$\mu(f|_X)=$dim$_\mathbb{C}\frac{\mathcal{O}_{\mathcal{X},0}}{\langle
s
\rangle+J(f_s,\phi_s)}=$dim$_\mathbb{C}\frac{\frac{\mathcal{O}_{n+1}}{\langle
\Phi \rangle}}{\langle s \rangle+J(f_s,\phi_s)}=$
dim$_\mathbb{C}\frac{\mathcal{O}_n}{\langle \phi
\rangle+J(f,\phi)}=\mu(X,0)+\mu(X\cap f^{-1}(0),0)$, where
$\mathcal{O}_{k}$ denotes the local ring of the function germs from
$(\mathbb{C}^k,0)$ to $\mathbb{C}$ and the last equality is the
Lê-Greuel formula.

In order to produce a result like Theorem \ref{teorema 1.1 greuel}
for this Milnor number, we use the following result of Teissier,
which gives us a characterization for the integral closure of an
ideal of the ring $\mathcal{O}_{X,x}$, where $(X,x)$ is an analytic
variety.

\begin{teo} \cite[Proposition 0.4]{Teissier} \label{criterio avaliativo} Let $(X,x)\subset(\mathbb{C}^n,0)$ be an analytic variety and let $I$ be an ideal in $\mathcal{O}_{X,x}$, the following conditions are equivalent:

        (i) $h \in \overline{I}$, where $\overline{I}=\{h \in \mathcal{O}_{\mathbb{C}\times X} \ | \ \exists \ a_i \in I^i$ with $h^k+a_1h^{k-1}+\ldots+a_{k-1}h+a_k=0\}$.

        (ii) For each system of generators $h_1,\ldots,h_r$ of $I$ there exists
        a neighborhood $U$ of $0$ and a constant $c>0$ such that $|h(x)|\leq
        csup\{|h_1(x)|,\ldots,|h_r(x)|\},$ $\forall \ x \in U$.

        (iii) (Valuation Criterion) For each analytic curve $\gamma: (\mathbb{C},0)\rightarrow
        (X,x)$, $h\circ \gamma \in
        (\gamma^*(I))\mathcal{O}_1$, where $(\gamma^*(I))\mathcal{O}_1$ is the
        ideal generated by $h_i\circ \gamma$, $i=1,\ldots,r$.

         $(iv)$ $\nu(h\circ\gamma)\geq
    inf\{\nu(h_1\circ\gamma),\ldots,\nu(h_r\circ\gamma)\}$, for $\nu$
    being the usual valuation of a complex curve.

\end{teo}

\section{Main results}


Let $(X,0)\subset (\mathbb{C}^n,0)$ be the ICIS defined by a
holomorphic map germ $\phi:(\mathbb{C}^n,0)\rightarrow
    (\mathbb{C}^p,0)$ and let $f:(X,0)\rightarrow (\mathbb{C},0)$ be
    a holomorphic function germ with isolated singularity. Consider

\begin{center} $\begin{array}{cccc}
  F: & (\mathbb{C}\times X,0) & \rightarrow & (\mathbb{C},0) \\
   & (t,x) & \mapsto & f_t(x)
\end{array}$ \end{center} a (flat) deformation of $f$. We say that $F$ is $\mu$-constant if $\mu(f_t|_X)=\mu(f|_X)$ for $t$
small enough.

The main goal of this work is to study when $F$ is $\mu$-constant,
that is, to present a result like Theorem \ref{teorema 1.1 greuel}
for a family of functions on an ICIS.

The assertions of the Theorem \ref{teorema 1.1 greuel} in this
context would be

(1$_X$) $F$ is $\mu$-constant.

   (2$_X$)  For every holomorphic curve $\gamma: (\mathbb{C},0)\rightarrow
   (\mathbb{C}\times X,0)$

   \begin{center} $\nu\left((\partial F/\partial
           t)\circ\gamma\right)>\inf\{\nu(B_i\circ\gamma) \ | \ i=1,\ldots,r\}$ in
       $\mathcal{O}_{\mathbb{C}\times X}$, \end{center}

       (where $\nu$ denotes the usual valuation of a complex curve).

   (3$_X$) Same statement as in (2) with $``>"$ replaced by $``\geq"$.

   (4$_X$) $\partial F/\partial t \in \overline{J_X}$
   in $\mathcal{O}_{\mathbb{C}\times X}$.

   (5$_X$) $\partial F/\partial t \in \sqrt{J_X}$ in
   $\mathcal{O}_{\mathbb{C}\times X}$.

    (6$_X$) $v(J_X)=\{(t,x) \in \mathbb{C}\times\mathbb{C}^n \ | \ B_i(t,x)=0, \ i=1,\ldots,r\}=\mathbb{C}\times\{0\}$ near $(0,0)$,
    where $B_1,\ldots,B_r$ are the maximal order minors of the Jacobian matrix of $(f_t,\phi)$ (partial
derivatives with respect to $x$, only) and $J_X$
    is the ideal generated by them.

\vspace{0.2cm}

Unfortunately, in this singular context, these assertions are not
equivalent. But in this section we show

\begin{center} (2$_X$) $\Rightarrow$ (3$_X$) $\Leftrightarrow$ (4$_X$)
$\Rightarrow$ (5$_X$),

(1$_X$) $\Leftrightarrow$ (6$_X$),

(6$_X$) $\Rightarrow$ (5$_X$),

(4$_X$) $\Rightarrow$ (1$_X$). \end{center}

Also, we present examples for

\begin{center}

(1$_X$) $\nRightarrow$ (2$_X$),

(1$_X$) $\nRightarrow$ (3$_X$),

(3$_X$) $\nRightarrow$ (2$_X$),

(5$_X$) $\nRightarrow$ (2$_X$),

(5$_X$) $\nRightarrow$ (1$_X$),

(5$_X$) $\nRightarrow$ (3$_X$).

\end{center}

Although the integral closure that appears in the assertion (4$_X$)
is the one of an ideal we need to work with integral closure of
submodule to show that (4$_X$) implies (1$_X$). Therefore, we
remember, now, its definition.

\begin{defi} \cite[Definition 1.3]{Gaffney} Suppose $(X,0)$ is a complex analytic germ, $M$ a submodule of $\mathcal{O}_{X,0}^p$. Then $h \in
\mathcal{O}_{X,0}^p$ is in the integral closure of $M$, denoted by
$\overline{M}$, if and only if for all
$\gamma:(\mathbb{C},0)\rightarrow(X,0)$, $h\circ \gamma \in
(\gamma^*(M))\mathcal{O}_1$.
\end{defi}

Replacing $\mathcal{O}_1$ by its maximal ideal $\mathcal{M}_1$ we
get the definition of strict integral closure of $M$, which is
denoted by $\overline{M}^{\dagger}$ (see \cite[Definition
1.1]{Gaffney Equisingularity}). In this case, $h\in
\overline{M}^{\dagger}$ it is said strictly dependent on $M$.

\vspace{0.2cm}

Our next result provides a way to ensure the constancy of $\mu(X\cap
f_t^{-1}(0),0)$.

 \begin{teo} \label{lema 4 implica 1} Let $(X,0)\subset
(\mathbb{C}^n,0)$ be an ICIS defined by
$\phi:(\mathbb{C}^n,0)\rightarrow
    (\mathbb{C}^p,0)$ and let $f:(X,0)\rightarrow (\mathbb{C},0)$ be a holomorphic germ with isolated singularity. Consider $F:(\mathbb{C}\times X,0)\rightarrow
    (\mathbb{C},0)$ a deformation of $f$ and $G:\mathbb{C}\times\mathbb{C}^n\rightarrow
\mathbb{C}\times\mathbb{C}^p$ given by
    $G(t,x)=(F(t,x),\phi(x))$. If $\partial F/\partial t \in \overline{J_X}$ in $\mathcal{O}_{\mathbb{C}\times X}$
    then:

(i) $\left(\partial F/\partial t,0\right) \in \overline{\{
x_i(\partial G/\partial x_j)\}}$ in $\mathcal{O}_{X\cap
f_t^{-1}(0)}$ with $i,j=1,\ldots,n$, where $\overline{\{
x_i(\partial G/\partial x_j)\}}$ denotes the integral closure of the
$\mathcal{O}_X$-module $\{ x_i(\partial G/\partial x_j)\}$.

(ii) $X\cap f_t^{-1}(0)$ is an ICIS with $\mu(X\cap f_t^{-1}(0),0)$
constant.

 \end{teo}

\noindent \textbf{Proof:} (i) We remember that $J_X$ is the ideal in
$\mathcal{O}_{\mathbb{C}\times X}$ generated by the minors of order
$p+1$ of the matrix

\begin{center} $\left(
                                               \begin{array}{cccc}
                                                 {f_t}_{x_1} & {f_t}_{x_2} & \cdots & {f_t}_{x_n} \\
                                                 {\phi_1}_{x_1} & {\phi_1}_{x_2} & \cdots & {\phi_1}_{x_n} \\
                                                 \vdots & \vdots & \vdots & \vdots \\
                                                 {\phi_p}_{x_1} & {\phi_p}_{x_2} & \cdots & {\phi_p}_{x_n} \\
                                               \end{array}
                                             \right)$, where $\phi=(\phi_1,\ldots,\phi_p)$ and $f_t(x)=F(t,x)$. \end{center}

These minors are

\begin{center} $B_v=\left|
                                               \begin{array}{cccc}
                                                 {f_t}_{x_{j_1}} & {f_t}_{x_{j_2}} & \cdots & {f_t}_{x_{{j_{p+1}}}} \\
                                                 {\phi_1}_{x_{j_1}} & {\phi_1}_{x_{j_2}} & \cdots & {\phi_1}_{x_{j_{p+1}}} \\
                                                 \vdots & \vdots & \vdots & \vdots \\
                                                 {\phi_p}_{x_{j_1}} & {\phi_p}_{x_{j_2}} & \cdots & {\phi_p}_{x_{j_{p+1}}} \\
                                               \end{array}
                                             \right|={f_t}_{x_{j_1}}\left|
                                               \begin{array}{cccc}
                                                 {\phi_1}_{x_{j_2}} & {\phi_1}_{x_{j_3}} & \cdots & {\phi_1}_{x_{j_{p+1}}} \\
                                                 \vdots & \vdots & \vdots & \vdots \\
                                                 {\phi_p}_{x_{j_2}} & {\phi_p}_{x_{j_3}} & \cdots & {\phi_p}_{x_{j_{p+1}}} \\
                                               \end{array}
                                             \right|-{f_t}_{x_{j_2}}\left|
                                               \begin{array}{cccc}
                                                 {\phi_1}_{x_{j_1}} & {\phi_1}_{x_{j_3}} & \cdots & {\phi_1}_{x_{j_{p+1}}} \\
                                                 \vdots & \vdots & \vdots & \vdots \\
                                                 {\phi_p}_{x_{j_1}} & {\phi_p}_{x_{j_3}} & \cdots & {\phi_p}_{x_{j_{p+1}}} \\
                                               \end{array}
                                             \right|+\ldots+(-1)^{p+2}{f_t}_{x_{j_{p+1}}}\left|
                                               \begin{array}{cccc}
                                                 {\phi_1}_{x_{j_1}} & {\phi_1}_{x_{j_2}} & \cdots & {\phi_1}_{x_{j_{p}}} \\
                                                 \vdots & \vdots & \vdots & \vdots \\
                                                 {\phi_p}_{x_{j_1}} & {\phi_p}_{x_{j_2}} & \cdots & {\phi_p}_{x_{j_{p}}} \\
                                               \end{array}
                                             \right|$,  \end{center}
                                             for each vector $v=(j_1,\ldots,j_{p+1})$ with $j_1<\ldots<j_{p+1}$ and $j_1,\ldots,j_{p+1}\in\{1,\ldots,n\}$.

Let $A_{v_i}=\left|
                                      \begin{array}{cccccc}
                                         {\phi_1}_{x_{j_1}} & \cdots &  {\phi_1}_{x_{{j_{i-1}}}} &  {\phi_1}_{x_{j_{i+1}}} &  \cdots &  {\phi_1}_{x_{j_{p+1}}} \\
                                        \vdots & \vdots & \vdots & \vdots & \vdots & \vdots \\
                                        {\phi_p}_{x_{j_1}} & \cdots &  {\phi_p}_{x_{j_{i-1}}} &  {\phi_p}_{x_{j_{i+1}}} &  \cdots &  {\phi_p}_{x_{j_{p+1}}} \\
                                      \end{array}
                                    \right|$, for each $i=1,\ldots,p+1$.

Thus
$B_v={f_t}_{x_{j_1}}A_{v_1}-{f_t}_{x_{j_2}}A_{v_2}+\ldots+(-1)^{1+p+1}{f_t}_{x_{j_{p+1}}}A_{{v_{p+1}}}$.
Then

\noindent$\displaystyle\left(\displaystyle\sum_{i=1}^{p+1}(-1)^{1+i}{f_t}_{x_{j_i}}{A_{v_i}},\displaystyle\sum_{i=1}^{p+1}(-1)^{1+i}{\phi_1}_{x_{j_i}}{A_{v_i}},\displaystyle\sum_{i=1}^{p+1}(-1)^{1+i}{\phi_2}_{x_{j_i}}{A_{v_i}},\ldots,\displaystyle
\displaystyle\sum_{i=1}^{p+1}(-1)^{1+i}{\phi_p}_{x_{j_i}}{A_{v_i}}\right)
=(B_v,0,\ldots,0)$, because all the entries, but the first, is the
determinant of a matrix with two equal lines.

Thus, $(B_v,0,\ldots,0)\in\mathcal{M}_nL_i$ for $i=1,\ldots,n$,
where $L_i:=\partial G/\partial
x_i=({f_t}_{x_i},{\phi_1}_{x_i},\ldots,{\phi_p}_{x_i})$ and
$\mathcal{M}_n=\{h \in \mathcal{O}_n \ | \ h(0)=0\}$.

By the hypothesis, $\partial F/\partial t \in \overline{J_X}$ in
$\mathcal{O}_{\mathbb{C}\times X}$ then, by Theorem \ref{criterio
avaliativo}, for all $\varphi:(\mathbb{C},0)\rightarrow
\mathbb{C}\times X$, $(\partial F/\partial t)\circ\varphi \in
\langle B_v\circ\varphi\rangle$ and, therefore, $\left((\partial
F/\partial t)\circ\varphi,0\right)\in\langle
(B_v\circ\varphi,0)\rangle\subset(\mathcal{M}_nL_i)\circ\varphi$,
for $v=(j_1,\ldots,j_{p+1})$ with $j_1<\ldots<j_{p+1}$ and
$j_1,\ldots,j_{p+1}\in\{1,\ldots,n\}$. Thus, $\left(\partial
F/\partial t,0\right)\in \overline{\{ x_i(\partial G/\partial
x_j)\}}$.

\vspace{0.3cm}

(ii) By the item (i), $\left(\partial F/\partial t,0\right) \in
\overline{\{x_i(\partial G/\partial x_j)\}}$ in $\mathcal{O}_{X\cap
f_t^{-1}(0)}$. Then, $\{X\cap f_t^{-1}(0)\}$ is Whitney regular (see
\cite[Theorem 2.5]{Gaffney}). In particular, $\mu(X\cap
f_t^{-1}(0),0)$ is constant (see \cite[Theorem 6.1]{Gaffney 1}).
\fim

\vspace{0.3cm}

From the proof of the item (ii) of the previous theorem,
$\left(\partial F/\partial t,0\right) \in \overline{\{x_i(\partial
G/\partial x_j)\}}$ implies that $\mu(X\cap f_t^{-1}(0),0)$ is
constant.

We see, in the following example that $\left(\partial F/\partial
t,0\right) \in \overline{\{\partial G/\partial x_j\}}$ does not imply
that $F$ is $\mu$-constant

\begin{ex} \label{contra exemplo do nosso teorema} Let $(X,0)\subset
(\mathbb{C}^3,0)$ be the zero set of $\phi(x,y,z)=x^5+y^3+z^2$ and
let $f:(X,0)\rightarrow(\mathbb{C},0)$ be defined by $f(x,y,z)=xy$.
Consider the deformation of $f$, $F(t,(x,y,z))=xy-tz$. Note that
$J_X=\langle -2xz-3ty^2,-2yz-5tx^4,3y^3-5x^5 \rangle$.

Now let's consider $\gamma:(\mathbb{C},0)\rightarrow
(\mathbb{C}\times X,0)$ given by $\gamma(s)=(0,-s^2,0,s^5)$. We have
that $\nu \left((\partial F/\partial t)\circ \gamma\right)=5$ and
$\nu(J_X\circ\gamma)=7$.

Therefore, $\nu \left((\partial F/\partial t)\circ
\gamma\right)<\nu(J_X\circ\gamma).$ Thus, $\partial F/\partial t
\notin \overline{J_X}$.

On the other hand, we have that

\begin{center} $\left(\partial F/\partial t,0\right)=(-z,0)$ and $\{\partial
G/\partial x,\partial G/\partial y,\partial G/\partial
z\}=\{(y,5x^4),(x,3y^2),(-t,2z)\}$, \end{center} where
$G(t,(x,y,z))=(F(t,(x,y,z)),\phi(x,y,z))$.

Note that

\begin{center} $(-z,0) \in\overline{\{\partial G/\partial x,\partial
G/\partial y,\partial G/\partial z\}}$ if and only if
$J_2((-z,0),\{\partial G/\partial x,\partial G/\partial y,\partial
G/\partial z\})\subset \overline{J_2(\{\partial G/\partial
x,\partial G/\partial y,\partial G/\partial z\})}$ \end{center} (see
\cite[Proposition 1.7]{Gaffney}), where $J_2((-z,0),\{\partial
G/\partial x,\partial G/\partial y,\partial G/\partial z\})$ are the
minors $2\times 2$ of matrix

\begin{center} $\left(
             \begin{array}{cccc}
               -z & y & x & -t\\
               0 & 5x^4 & 3y^2 & 2z\\
             \end{array}
           \right)
$.\end{center}

Therefore, we just show that

\begin{center} $-5x^4z, -3y^2z$ and $-2z^2$ belong to $\overline{\langle
-2xz-3ty^2,-2yz-5tx^4,3y^3-5x^5 \rangle}$. \end{center}

Let $\gamma:(\mathbb{C},0)\rightarrow(\mathbb{C}\times X,0)$,
$\gamma(s)=(\gamma_1(s),\gamma_2(s),\gamma_3(s),\gamma_4(s))$. Thus,
$\gamma_2^5+\gamma_3^3+\gamma_4^2=0$.

By the Theorem \ref{criterio avaliativo} we have the desired.

Therefore, $\left(\partial F/\partial t,0\right) \in
\overline{\{\partial G/\partial x,\partial G/\partial y,\partial
G/\partial z\}}$.

However, $\mu(f_0|_X)=17$ and $\mu(f_t|_X)=16$, for $t\neq 0$.

 \end{ex}

We are now ready to prove our main result.

\begin{teo} \label{teo demonstracao} Let $(X,0)\subset
(\mathbb{C}^n,0)$ be an ICIS defined by
$\phi:(\mathbb{C}^n,0)\rightarrow
    (\mathbb{C}^p,0)$, let $f:(X,0)\rightarrow (\mathbb{C},0)$ be
    a germ with isolated singularity and let $f_t:(X,0)\rightarrow
    (\mathbb{C},0)$ be a (flat) deformation of $f$. Then:

\begin{center} \rm{(2$_X$) $\Rightarrow$ (3$_X$) $\Leftrightarrow$ (4$_X$)
$\Rightarrow$ (5$_X$),

(1$_X$) $\Leftrightarrow$ (6$_X$),

(6$_X$) $\Rightarrow$ (5$_X$),

(4$_X$) $\Rightarrow$ (1$_X$)}. \end{center}

\end{teo}

\noindent \textbf{Proof:} (2$_X$) $\Rightarrow$ (3$_X$) is trivial,
(3$_X$) $\Leftrightarrow$ (4$_X$) follows from $(i) \Leftrightarrow
(iv)$ in Theorem \ref{criterio avaliativo} and (4$_X$) $\Rightarrow$
(5$_X$) because $\overline{J_X}\subset \sqrt{J_X}$.

\vspace{0.3cm}

(1$_X$) $\Rightarrow$ (6$_X$):

Since $\frac{\mathcal{O}_n}{\langle \phi\rangle+J_X}$ is
Cohen-Macaulay, by the principle of conservation of number and by
the hypothesis, we have that
\begin{center} $\mu(f_t|_X,0)=\mu(f|_X,0)=\displaystyle\sum_{(t,x)\in
\{t\}\times v(J_X)}\mu({f_t}|_X,x)$, \end{center} then
$\mu(f_t|_X,x)=0$ for all $x\neq 0$ and therefore
$v(J_X)=\mathbb{C}\times\{0\}$ near $(0,0)$.

\vspace{0.3cm}

(6$_X$) $\Rightarrow$ (1$_X$):

By the hypothesis $v(J_X)=\mathbb{C}\times\{0\}$ near $(0,0)$.
Again, by the principle of conservation of number,
$\mu(f|_X,0)=\mu(f_t|_X,0)$. Therefore, $F$ is $\mu$-constant.

\vspace{0.3cm}

(6$_X$)$\Rightarrow$ (5$_X$):

By the hypothesis $v(J_X)=\mathbb{C}\times\{0\}$ near $(0,0)$. We
have that $F\in \mathcal{M}_n$, then ${(\partial F/\partial
t)}|_{v(J_X)}={(\partial F/\partial
t)}|_{\mathbb{C}\times\{0\}}\equiv 0$. Thus, $v(J_X)\subset
\displaystyle v\left(\partial F/\partial t\right)$. Therefore,
$\partial F/\partial t\in \sqrt{J_X}$ by the Hilbert Nullstellensatz
Theorem.

\vspace{0.3cm}

(4$_X$) $\Rightarrow$ (1$_X$):

By the Lê-Greuel formula, $\mu(f_t|_X)=\mu(X,0)+\mu(X\cap
f_t^{-1}(0),0)$. Therefore, since, $\mu(X\cap f_t^{-1}(0),0)$ is
constant by the Theorem \ref{lema 4 implica 1}, $\mu(f_t|_X)$ is
also constant. \fim

\vspace{0.2cm}

From here to the end of this section, our goal is to present
counterexamples for the other implications.

A good idea to look for such counterexamples is to see families
which we know, at first, that have constant Milnor number. For this,
we recall the results of \cite{Bruna} on deformations of weighted
homogeneous germs.

\vspace{0.3cm}

We say that
$f=(f_1,\ldots,f_p):(\mathbb{C}^n,0)\rightarrow(\mathbb{C}^p,0)$ is
weighted homogeneous of type $(w_1,\ldots,w_n;d_1,\ldots,d_p)$ if
given $(w_1,\ldots,w_n;d_1,\ldots,d_p)$ with $w_i, d_j\in
\mathbb{Q}^+$ we have that for all $\lambda \in \mathbb{C}-\{0\}$:

\begin{center} $f(\lambda^{w_1}x_1,\ldots,\lambda^{w_n}x_n)=(\lambda^{d_1}f_1(x),\ldots,\lambda^{d_p}f_p(x))$. \end{center}

We call $d_j$ the weighted degree of $f_j$, which is denoted by
$wt(f_j)$ and $w_i$ is called weight of the variable $x_i$.

Let $(X,0)\subset (\mathbb{C}^n,0)$ be a germ of analytic variety
defined by $\phi:(\mathbb{C}^n,0)\rightarrow (\mathbb{C}^p,0)$. If
$\phi$ is weighted homogeneous we say $(X,0)$ is weighted
homogeneous.

A non-zero polynomial germ $f:(\mathbb{C}^n,0)\rightarrow\mathbb{C}$
can be written

\begin{center} $f=f_d+f_{d+1}+\ldots+f_l$, \end{center} where $f_d\neq 0$ and each $f_i$ is weighted homogeneous of degree
$i$. We say that $f_d$ is the initial part of $f$, which is
denoted by $in(f)$.

Moreover, if $f:(\mathbb{C}^n,0)\rightarrow(\mathbb{C},0)$ is a
polynomial function and
$f_t:(\mathbb{C}^n,0)\rightarrow(\mathbb{C},0)$ a deformation of $f$
we have that $f_t$ can be written as:

\begin{center} $f_t(x)=f(x)+\displaystyle\sum_{i=1}^k\sigma_i(t)\alpha_i(x),$ \end{center} with $\alpha_i:(\mathbb{C}^n,0)\rightarrow (\mathbb{C},0)$ and
$\sigma_i:(\mathbb{C},0)\rightarrow (\mathbb{C},0)$.

If $wt(in(\alpha_i))\geq wt(in(f))$ for $i=1,\ldots,k$ we say that
$f_t$ is a non-negative deformation of $f$.

\begin{teo} \cite[Theorem 4.4]{Bruna} \label{deformacao acima} Let $(X,0)\subset(\mathbb{C}^n,0)$ be a weighted homogeneous ICIS and let $f:(\mathbb{C}^n,0)\rightarrow \mathbb{C}$
be a weighted homogeneous germ with an isolated singularity with the
same weights of $(X,0)$ and let $f_t$ be a deformation of $f$. If
$f_t$ is non-negative, then $\mu(f_t|_X)$ is constant.
\end{teo}

Using this result we show in the following example that (1$_X$) does
not imply (2$_X$) and (1$_X$) does not imply (3$_X$).

\begin{ex} Let $(X,0)\subset (\mathbb{C}^2,0)$ be defined by $\phi(x,y)=x^p-y^q$, with $q\geq 3$ and
$p>q$, and let $f:(X,0)\rightarrow(\mathbb{C},0)$ be defined by
$f(x,y)=x$. Consider the deformation of $f$ defined by
$F(t,(x,y))=x+ty$. Note that $J_X=\langle -qy^{q-1}-ptx^{p-1}
\rangle$.

Now let's consider $\gamma:(\mathbb{C},0)\rightarrow
(\mathbb{C}\times X,0)$ given by $\gamma(s)=(0,s^q,s^p)$ and we have
that $\nu \left((\partial F/\partial t)\circ \gamma\right)=p$ and
$\nu(J_X\circ\gamma)=(q-1)p$.

Therefore, $\nu \left((\partial F/\partial t)\circ
\gamma\right)<\nu(J_X\circ\gamma).$ Thus, (2$_X$) and (3$_X$) are
not true.

On the other hand since $f_t$ is a non-negative deformation,
$\mu(f_t|_X)$ is constant by the Theorem \ref{deformacao acima}.
Thus, (1$_X$) is true.

\end{ex}

We show in the next example that (3$_X$) does not imply (2$_X$) and
(5$_X$) does not imply (2$_X$).

\begin{ex} Let $(X,0)\subset (\mathbb{C}^2,0)$ be defined by $\phi(x,y)=x^{2q}-y^q$ with $q\geq 2$ and let $f:(X,0)\rightarrow(\mathbb{C},0)$ be
defined by $f(x,y)=x^{2q}+y^q$. Consider the deformation of $f$
given by $F(t,(x,y))=x^{2q}+y^q+tx^{4q-3}$. Note that $J_X=\langle
-4q^2x^{2q-1}y^{q-1}-q(4q-3)tx^{4q-4}y^{q-1} \rangle$.

Let $\gamma:(\mathbb{C},0)\rightarrow(\mathbb{C}\times X,0)$,
$\gamma(s)=(\gamma_1(s),\gamma_2(s),\gamma_3(s))$ such that
$\gamma_2^{2q}-\gamma_3^q=0$. Thus, $\gamma_2^{2q}=\gamma_3^q$.
Then, $2\nu(\gamma_2)=\nu(\gamma_3)$.

Since $\partial F/\partial t=x^{4q-3}$, we have that $\nu
\left((\partial F/\partial t)\circ
\gamma\right)=(4q-3)\nu(\gamma_2)$. Furthermore
$\nu(J_X\circ\gamma)=\nu(-4q^2
\gamma_2^{2q-1}\gamma_3^{q-1}-q(4q-3)\gamma_1\gamma_2^{4q-4}\gamma_3^{q-1})=(2q-1)\nu(\gamma_2)+2(q-1)\nu(\gamma_2)=(4q-3)\nu(\gamma_2)$.

Thus, (3$_X$) is true, consequently (4$_X$) is true and (5$_X$) is
true.

Now, let's consider $\gamma:(\mathbb{C},0)\rightarrow
(\mathbb{C}\times X,0)$ given by $\gamma(s)=(0,s,s^2)$ and we have
that $\nu \left((\partial F/\partial t)\circ \gamma\right)=4q-3$ and
$\nu(J_X\circ\gamma)=4q-3$. Thus, (2$_X$) is not true.

\end{ex}

The following example shows that (5$_X$) does not imply (1$_X$) and
(5$_X$) does not imply (3$_X$).

\begin{ex} Let $(X,0)\subset (\mathbb{C}^2,0)$ be defined by $\phi(x,y)=x^p-y^q$, with $q\geq 3$, $p<q$ and let
$f:(X,0)\rightarrow(\mathbb{C},0)$ be defined by $f(x,y)=x$.
Consider the deformation of $f$ given by $F(t,(x,y))=x+ty$. Note
that $J_X=\langle -qy^{q-1}-ptx^{p-1} \rangle$.

We have that $\mu(f|_X)=pq-p$ and $\mu(f_t|_X)=pq-q$. Therefore,
(1$_X$) is not true.

Now, let's consider $\gamma:(\mathbb{C},0)\rightarrow
(\mathbb{C}\times X,0)$ given by $\gamma(s)=(0,s^q,s^p)$. Thus $\nu
\left((\partial F/\partial t) \circ \gamma\right)=p$ and
$\nu(J_X\circ\gamma)=(q-1)p$. Therefore, $\nu \left((\partial
F/\partial t)\circ \gamma\right)<\nu(J_X\circ\gamma).$

On the other hand, $\partial F/\partial t \in \sqrt{J_X}$ in
$\mathcal{O}_{\mathbb{C}\times X}$.

\end{ex}

\section{Some examples}

In this section, we apply our results to produce examples of
families of functions in an ICIS which have constant Milnor number
although they do not satisfy the hypotheses of the Theorem
\ref{deformacao acima}. For this we resort to the results about of
the Newton polyhedron (see \cite{Saia}).

Let $g \in \mathcal{O}_n$ then $g$ can be written as $g(x)=\sum
a_kx^k$. We define the support of $g$ as supp $g:=\{k \in
\mathbb{Z}^n \ | \ a_k\neq 0\}$ and for $I$ an ideal in
$\mathcal{O}_n$, we define supp $I:=\bigcup \{$supp $g \ | \ g \in
I\}$.

The convex hull in $\mathbb{R}_+^n$ of the set $\bigcup \{k+v \ | \
k \in \ $supp $I, v \in \mathbb{R}_+^n\}$ is called Newton
polyhedron of $I$ and denoted by $\Gamma_+(I)$. The union of all
compact faces of $\Gamma_+(I)$ we denote by $\Gamma(I)$.

Let $\Delta\subset \Gamma_+(I)$ be a finite subset, we define
$g_{\Delta}=\sum_{k \in \Delta}a_kx^k$ for any germ $g(x)=\sum
a_kx^k$.

Given $\Delta$ a face of $\Gamma_+(I)$, we denote $C(\Delta)$ the
cone of half-rays emanating from $0$ and passing through $\Delta$.
We define $C[[\Delta]]$, the ring of power series with non-zero
monomials $x^k=x_1^{k_1}x_2^{k_2}\cdot\ldots\cdot x_n^{k_n}$ such
that $k=(k_1,\ldots,k_n)\in C(\Delta)$. When the ideal generated by
${g_1}_{\Delta},{g_2}_{\Delta},\ldots,{g_n}_{\Delta}$ has finite
codimension in $C[[\Delta]]$ we say that the compact face $\Delta
\subset \Gamma(I)$ is Newton non-degenerate. Furthermore, if all
compact faces of $\Gamma(I)$ are Newton non-degenerate then the
ideal $I$ is said Newton non-degenerate.

Equivalently, in \cite[p.2]{Saia}: $I$ is Newton non-degenerate if
for each compact face $\Delta\subset\Gamma(I)$, the equations
${g_1}_{\Delta}(x)={g_2}_{\Delta}(x)=\ldots={g_n}_{\Delta}(x)=0$
have no common solution in $(\mathbb{C}-\{0\})^n$.

\begin{teo} \cite[Theorem 3.4]{Saia} \label{fecho integral e poliedro de newton} Let $I=\langle g_1,g_2,\ldots,g_s \rangle$ be an
ideal of finite codimension in $\mathcal{O}_n$. Then $I$ is Newton
non-degenerate if and only if $\Gamma_+(I)=C(\overline{I})$, where
$C(\overline{I})$ is the convex hull in $\mathbb{R}_+^n$ of the set
$\bigcup \{m \ | \ x^m \in \overline{I}\}$.
\end{teo}

Let $(X,0)\subset(\mathbb{C}^n,0)$ be the ICIS defined by the zero
set of a map germ
$\phi:(\mathbb{C}^n,0)\rightarrow(\mathbb{C}^p,0)$. Let
$f:(X,0)\rightarrow (\mathbb{C},0)$ be a holomorphic function germ
with isolated singularity and $F:(\mathbb{C}\times X,0)\rightarrow
(\mathbb{C},0)$ a deformation of $f$ given by $F(t,x)=f(x)+tg(x)$,
where $g$ is a holomorphic function germ such that $g(0)=0$.

We denote by $B_1,\ldots,B_r$ the minors of order $p+1$ of
$J(f_t,\phi)$, $B_1^0,\ldots,B_r^0$ the minors of order $p+1$ of
$J(f,\phi)$ and $C_1,\ldots,C_r$ the minors of order $p+1$ of
$J(g,\phi)$. Then, $B_i=B_i^0+tC_i$, by the multilinearity of the
determinant, for $i=1,\ldots,r$.

With the above notation,

\begin{lem} \label{estar no fecho} If $g\in \overline{\langle\phi\rangle+J(f,\phi)}$ in $\mathcal{O}_{n+1}$ and $C_i\in\overline{\langle\phi\rangle+J(f,\phi)}$ in
$\mathcal{O}_{n+1}$ for $i=1,\ldots,r$ then $g\in \overline{J_X}$ in
$\mathcal{O}_{\mathbb{C}\times X}$. \end{lem}

\noindent \textbf{Proof:} Since
$C_i\in\overline{\langle\phi\rangle+J(f,\phi)}$ then, by Theorem
\ref{criterio avaliativo}, there exists a neighborhood $U$ of $0$
and a constant $c>0$ such that

\begin{center} $|t||C_i|\leq |t|csup_i\{|\phi|,|B_i^0|\}$. \end{center}

Besides that,

\begin{center}
$sup_i\{|\phi|,|B_i|\}=sup_i\{|\phi|,|B_i^0+tC_i|\}\geq$

\vspace{0.2cm}

$sup_i\{|\phi|,|B_i^0|\}-|t|sup_i\{|\phi|,|C_i|\}\geq$

\vspace{0.2cm}

$sup_i\{|\phi|,|B_i^0|\}-|t|csup_i\{|\phi|,|B_i^0|\}\geq$

\vspace{0.2cm}

$(1-\alpha)sup_i\{|\phi|,|B_i^0|\}$, \end{center} where $0<\alpha<1$
and $|t|\leq \frac{\alpha}{c}$.

Thus, $sup_i\{|\phi|,|B_i|\}\geq Ksup_i\{|\phi|,|B_i^0|\}$ for $K>0$
and $t$ small enough. Therefore,
$\overline{\langle\phi\rangle+J(f,\phi)}\subseteq
\overline{\langle\phi\rangle+J_X}$ in $\mathcal{O}_{n+1}$. Then
$g\in \overline{\langle\phi\rangle+J_X}$ in $\mathcal{O}_{n+1}$.
Thus, by Theorem \ref{criterio avaliativo} for all
$\gamma:\mathbb{C}\rightarrow \mathbb{C}\times\mathbb{C}^n$ we have
that $g\circ\gamma \in \langle \phi\circ\gamma, B_1\circ
\gamma,\ldots,B_r\circ \gamma \rangle$. Then, for all
$\gamma:\mathbb{C}\rightarrow \mathbb{C}\times
X\subset\mathbb{C}\times\mathbb{C}^n$ we have that $g\circ\gamma \in
\langle B_1\circ \gamma,\ldots,B_r\circ \gamma \rangle$. Again by
Theorem \ref{criterio avaliativo} we concluded that $g\in
\overline{J_X}$ in $\mathcal{O}_{\mathbb{C}\times X}$. \fim

\vspace{0.2cm}

With this lemma, in the following corollaries, we relate the Newton
polyhedron to the item (4$_X$) of our Theorem \ref{teo
demonstracao}.

\begin{cor} \label{estar no fecho 1} If $\Gamma_+(g)\subset C(\overline{\langle\phi\rangle+J(f,\phi)})$ in $\mathcal{O}_{n+1}$ and
$\Gamma_+(C_i)\subset C(\overline{\langle\phi\rangle+J(f,\phi)})$ in
$\mathcal{O}_{n+1}$ for $i=1,\ldots,r$ then $g\in \overline{J_X}$ in
$\mathcal{O}_{\mathbb{C}\times X}$. \end{cor}

\noindent \textbf{Proof:} Since $\Gamma_+(g)\subset
C(\overline{\langle\phi\rangle+J(f,\phi)})$ in $\mathcal{O}_{n+1}$
and $\Gamma_+(C_i)\subset
C(\overline{\langle\phi\rangle+J(f,\phi)})$ in $\mathcal{O}_{n+1}$
then $g \in \overline{\langle\phi\rangle+J(f,\phi)}$ in
$\mathcal{O}_{n+1}$ and $C_i\in
\overline{\langle\phi\rangle+J(f,\phi)}$ in $\mathcal{O}_{n+1}$.

Thus, by the Lemma \ref{estar no fecho}, $g\in \overline{J}_X$ in
$\mathcal{O}_{\mathbb{C}\times X}$. \fim

\begin{cor} \label{estar no fecho 2} If $\langle\phi\rangle+J(f,\phi)$ is Newton non-degenerate, $\Gamma_+(g)\subset \Gamma_+(\langle\phi\rangle+J(f,\phi))$ in
 $\mathcal{O}_{n+1}$ and $\Gamma_+(C_i)\subset \Gamma_+(\langle\phi\rangle+J(f,\phi))$ in
 $\mathcal{O}_{n+1}$ for $i=1,\ldots,r$
 then $g\in \overline{J_X}$ in $\mathcal{O}_{\mathbb{C}\times X}$. \end{cor}

\noindent \textbf{Proof:} Since $\langle\phi\rangle+J(f,\phi)$ is
Newton non-degenerate then by Theorem \ref{fecho integral e poliedro
de newton} we have that
$\Gamma_+(\langle\phi\rangle+J(f,\phi))=C(\overline{\langle\phi\rangle+J(f,\phi)})$.

Therefore, by the Corollary \ref{estar no fecho 1}, $g\in
\overline{J}_X$ in $\mathcal{O}_{\mathbb{C}\times X}$. \fim

\vspace{0.2cm}

We are ready, now, to produce interesting examples of families with
constant Milnor number. We highlight that the families in the two
following examples do not satisfy the hypothesis of the Theorem
\ref{deformacao acima}.

\begin{ex} \label{ex nao quase homogenea 2} Let $(X,0)\subset
(\mathbb{C}^3,0)$ be defined by $\phi(x,y,z)=(xy,x^{15}+y^{10}+z^6)$
and let $f:(X,0)\rightarrow(\mathbb{C},0)$ be defined by
$f(x,y,z)=x+z$. We have that $\phi$ and $f$ are weighted homogeneous
with different weights. Therefore, we can not use the Theorem
\ref{deformacao acima}. Let $F:(\mathbb{C}\times
X)\rightarrow(\mathbb{C},0)$ be the deformation of $f$ defined by
$F(t,(x,y,z))=x+z+txy$.

We have that $\langle\phi\rangle+J(f,\phi)=\langle
xy,x^{15}+y^{10}+z^6,6xz^5+10y^{10}-15x^{15} \rangle$ is Newton
non-degenerate, $\Gamma_+(g)\subset
\Gamma_+(\langle\phi\rangle+J(f,\phi))$ in $\mathcal{O}_{n+1}$ and
$\Gamma_+(C_1)\subset \Gamma_+(\langle\phi\rangle+J(f,\phi))$ in
$\mathcal{O}_{n+1}$, where $g(x,y,z)=xy$ and $C_1$ is the
determinant of $J(g,\phi)$. Then, by Corollary \ref{estar no fecho
2}, $\partial F/\partial t\in \overline{J_X}$ in
$\mathcal{O}_{\mathbb{C}\times X}$. Therefore, by Theorem \ref{teo
demonstracao}, $F$ is $\mu$-constant.

\end{ex}

\begin{ex} \label{ex nao quase homogenea 3} Let $(X,0)\subset (\mathbb{C}^3,0)$ be defined by $\phi(x,y,z)=x^3+y^3+z^4+xyz$ and let $f:(X,0)\rightarrow(\mathbb{C},0)$ be defined by $f(x,y,z)=xy+z^2$. Consider $F:(\mathbb{C}\times
X)\rightarrow(\mathbb{C},0)$ be the deformation of $f$ defined by
$F(t,(x,y,z))=f(x,y,z)+tg(x,y,z)$. Note that $\phi$ is not weighted
homogeneous.

We have that $\partial F/\partial t=g$ and
$\langle\phi\rangle+J(f,\phi)=\langle x^3+y^3+z^4+xyz,
-x^2y+6y^2z+2xz^2-4xz^3,-xy^2+6x^2z+2yz^2-4yz^3,-3x^3+3y^3 \rangle$
is Newton non-degenerate. In addition, except for $g=z^3$ any other
deformation with a degree greater or equal to three is such that
$\Gamma_+(g)\subset \Gamma_+(\langle\phi\rangle+J(f,\phi))$ in
$\mathcal{O}_{n+1}$ and $\Gamma_+(C_i)\subset
\Gamma_+(\langle\phi\rangle+J(f,\phi))$ in $\mathcal{O}_{n+1}$,
where $C_i$ are the minors of order $2$ of $J(g,\phi)$ for
$i=1,2,3$. Thus, by Corollary \ref{estar no fecho 2}, $\partial
F/\partial t\in \overline{J_X}$ in $\mathcal{O}_{\mathbb{C}\times
X}$. Therefore, $F$ is $\mu$-constant.

\end{ex}

\section{Other invariants}

Given an ICIS $(X,0)\subset(\mathbb{C}^n,0)$, we have two important
subgroups of the $\mathcal{R}$ group of diffeomorphisms from
$(\mathbb{C}^n,0)$ to $(\mathbb{C}^n,0)$: one is the group
$\mathcal{R}_X$ of the diffeomorphisms which preserves $(X,0)$, the
other is $\mathcal{R}(X)$, the group of diffeomorphisms of $X$. We
know that if the germs $f,g:(X,0)\rightarrow(\mathbb{C},0)$ are
$\mathcal{R}_X$-equivalent then they are
$\mathcal{R}(X)$-equivalent, but the converse is not true.

In the smooth case, we know that a germ
$f:(\mathbb{C}^n,0)\rightarrow(\mathbb{C},0)$ is finitely determined
if and only if $\mu(f)$ is finite. There exist a generalization of this result
for the $\mathcal{R}_X$-group:
$f:(\mathbb{C}^n,0)\rightarrow(\mathbb{C},0)$ is
$\mathcal{R}_X$-finitely determined if and only if $\mu_{BR}(X,f)$ is finite.
Here $\mu_{BR}(X,f)$ is the Bruce-Roberts number defined in
\cite{Bruce e Roberts} in the following by

\begin{center} $\mu_{BR}(X,f)=\dim_{\mathbb{C}}\displaystyle\frac{\mathcal{O}_n}{J_f(\Theta_X)}$,
\end{center} where $\Theta_X$ is the $\mathcal{O}_n$-module of vector
fields in $(\mathbb{C}^n,0)$ which are tangent to $(X,0)$ and
$J_f(\Theta_X)=\langle df(\xi) \ | \ \xi\in \Theta_X \rangle$.

Because of this, it is important to know when a family has constant
Bruce-Roberts number.

Another important number related to
$f:(X,0)\rightarrow(\mathbb{C},0)$ is the Euler obstruction,
$Eu_{f,X}(0)$. This number is very studied for instance in
\cite{Brasselet} and \cite{Grulha}.

In this section we use our results to conclude when the
Bruce-Roberts number or the Euler obstruction of a family
$f_t:(X,0)\rightarrow(\mathbb{C},0)$ is constant. Also, we present a
sufficient condition for such a family to be $C^0-\mathcal{R}_X$(or
$C^0-\mathcal{R}(X))$-trivial.

In \cite{Tomazella}, Ahmed, Ruas and Tomazella studied the constancy
of the Bruce-Roberts number of a family. In order to understand
this, we need to know the logarithmic characteristic variety.

Let $U\subset \mathbb{C}^n$ be a neighborhood of origin. Suppose
that $\Theta_X=\langle \xi_1,\ldots,\xi_p \rangle$ em $U$. We define

\begin{center} $LC_U(X):=\{(x,\delta) \in T^*_U\mathbb{C}^n \ | \ \delta(\xi_i(x))=0, i=1,\ldots,p\}$, \end{center} where $T^*_U\mathbb{C}^n$
is the restriction of the cotangent bundle of $\mathbb{C}^n$ to $U$.

The logarithmic characteristic variety of $X$, which we denote by
$LC(X)$, is defined as the germ of $LC_U(X)$ in $T^*_0\mathbb{C}^n$,
the cotangent space of $\mathbb{C}^n$ in $0$.

\begin{teo} \label{teo cidinha, tomazella, ahmed} \cite[Theorem 3.11]{Tomazella} Suppose that $X$ is a hypersurface with isolated singularity and LC(X)
is Cohen-Macaulay, $f:(\mathbb{C}^n,0)\rightarrow (\mathbb{C},0)$
    a germ with isolated singularity and $F:(\mathbb{C}\times \mathbb{C}^n,0)\rightarrow
    (\mathbb{C},0)$ a deformation of $f$. If
$\partial F/\partial t\in
\overline{J_F(\Theta_X)}$ then $F$ is a $\mu_{BR}$-constant
deformation of $f$.
\end{teo}

In addition, given $(X,0)\subset(\mathbb{C}^n,0)$ an ICIS, Grulha in
\cite{Grulha} proved the equivalence between the constancy of
$\mu(X\cap f_t^{-1}(0),0)$ and $Eu_{f,X}(0)$, by the Lê-Greuel
formula this show the importance of our Theorem \ref{teo
demonstracao}.

\begin{teo} \cite[Proposition 5.17]{Grulha} \label{proposicao da constancia da obstrucao de euler} Let $(X,0)\subset
(\mathbb{C}^n,0)$ be an ICIS, and
$F:(\mathbb{C}^n\times\mathbb{C}^r,0)\rightarrow \mathbb{C}$ a
family of functions with isolated singularity, where
$F(t,x)=f_t(x)$. Then, the following statements are equivalent:

(1) $Eu_{f_t,X}(0)$ is constant for the family.

(2) $\mu(X\cap f_t^{-1}(0),0)$ is constant for the family.
\end{teo}

Using these results and the Theorem \ref{teo demonstracao} we obtain
the following applications of our study:

\begin{teo} \label{aplicacoes} Let $(X,0)\subset (\mathbb{C}^n,0)$ be an ICIS, $f:(X,0)\rightarrow (\mathbb{C},0)$
    a germ with isolated singularity $\mathcal{R}_X$-finitely determined and $F:(\mathbb{C}\times X,0)\rightarrow
    (\mathbb{C},0)$ a deformation of $f$. If $\partial F/\partial t\in \overline{J_X}$ then:

(i) $\tilde{F}$ is $C^0-\mathcal{R}_X$-trivial (and therefore
$\tilde{F}$ is $C^0-\mathcal{R}(X)$-trivial), where
$\tilde{F}:(\mathbb{C}\times\mathbb{C}^n,0)\rightarrow(\mathbb{C},0)$
is such that $F=\tilde{F}|_X$.

(ii) $Eu_{\tilde{f_t},X}(0)$ is constant, where
$\tilde{f}_t:(\mathbb{C}^n,0)\rightarrow(\mathbb{C},0)$ is such that
$f_t=\tilde{f}_t|_X$.

(iii) If $(X,0)$ is a hypersurface with isolated singularity and
LC(X) Cohen-Macaulay then $\mu_{BR}(X,\tilde{f}_t)$ is constant.

(iv) If $(X,0)$ is a weighted homogeneous hypersurface with isolated
singularity then $\mu_{BR}(X,\tilde{f}_t)$ is constant and
$m(\tilde{f}_t)$ is constant, where $m(\tilde{f}_t)$ is the
multiplicity of $\tilde{f}_t$.

\end{teo}

\noindent \textbf{Proof:}

(i) In \cite[Theorem 4.3]{cidinha e tomazella} it is shown that if
$(X,0)$ is an ICIS and $\partial F/\partial t\in
\overline{J_F(\Theta_X)}$ then $\tilde{F}$ is
$C^0-\mathcal{R}_X$-trivial. Therefore, as $J_X\subseteq
J_F(\Theta_X)$ we have the desired. Thus, $\tilde{F}$ also is
$C^0-\mathcal{R}(X)$-trivial.

(ii) Follows directly from Theorem \ref{teo demonstracao} and
Theorem \ref{proposicao da constancia da obstrucao de euler}.

(iii) Just use $J_X\subseteq J_F(\Theta_X)$ and Theorem \ref{teo
cidinha, tomazella, ahmed}.

(iv) In \cite[Theorem 4.2]{Bruna 1} it is shown that if $(X,0)$ is a
weighted homogeneous hypersurface with isolated singularity, then
$LC(X)$ is Cohen-Macaulay. Then follows as a consequence of item
(iii) that $\mu_{BR}(X,\tilde{f}_t)$ is constant. Thus,
$m(\tilde{f}_t)$ is constant (see \cite[Theorem 4.3]{Tomazella}).
\fim

\vspace{0.3cm}

As an application of the item (ii) of the Theorem \ref{aplicacoes}
we have that the families of Examples \ref{ex nao quase homogenea 2}
and \ref{ex nao quase homogenea 3} have constant Euler obstruction.

\section{Deformation of the ICIS}

Our final goal is to answer a natural question: What if we deform
the ICIS instead of the function germ?

Let $\phi:(\mathbb{C}^n,0)\rightarrow
    (\mathbb{C}^p,0)$ be a holomorphic map germ and let $(X,0)\subset (\mathbb{C}^n,0)$ be
    the ICIS defined by the zero set of $\phi$.

Let
$\Phi:(\mathbb{C}\times\mathbb{C}^n,0)\rightarrow(\mathbb{C}^p,0)$
be a deformation of $\phi$, defined by $\Phi(t,x)=\phi_t(x)$, such
that $\phi_0=\phi$ and $(X_t,0):=(\phi_t^{-1}(0),0)$ is an ICIS for
$t$ small enough. We denote $\mathcal{X}=\Phi^{-1}(0)$. Let
$f:(\mathbb{C}^n,0)\rightarrow
    (\mathbb{C},0)$ be a
holomorphic function germ such that

\begin{center} $\begin{array}{cccc}
  f_t: & (X_t,0) & \rightarrow & (\mathbb{C},0) \\
   & x & \mapsto & f(x)
\end{array}$ \end{center} has isolated singularity for $t$ small enough. We study here the
constancy of the Milnor number $\mu(f_t|_{X_t})$.

Consider $G:\mathbb{C}\times\mathbb{C}^n\rightarrow
\mathbb{C}\times\mathbb{C}^p$ defined by
    $G(t,x)=(f(x),\phi_t(x))$.

For the case where we did not deform $X$, we saw, in the Example
\ref{contra exemplo do nosso teorema}, that $\partial G/\partial t
\in \overline{\{\partial G/\partial x_j\}}$ in $\mathcal{O}_{X\cap
f_t^{-1}(0)}$ does not imply that $f_t$ is $\mu$-constant. We know
that if $\partial G/\partial t \in \overline{\{x_i(\partial
G/\partial x_j)\}}$ in $\mathcal{O}_{X_t\cap f_t^{-1}(0)}$ then
$(\mathcal{X}_0,\mathbb{C}\times\{0\})$ is Whitney regular (see
\cite[Theorem 2.5]{Gaffney}), where
$\mathcal{X}_0=\mathcal{X}-\mathbb{C}\times\{0\}$. We also know that
if $\partial G/\partial t \in \overline{\{\partial G/\partial
x_j\}}^{\dagger}$ then $A_{f_t}$ holds for the pair
$(\mathcal{X}_0,\mathbb{C}\times\{0\})$ (see \cite[Lemma
5.1]{Gaffney Specialization}) (see definition of $A_{f_t}$ in
\cite{Gaffney e Massey}). Hence, we can think about the relationship
between $\partial G/\partial t \in \overline{\{\partial G/\partial
x_j\}}^{\dagger}$ and the constancy of $\mu(f_t)$.

\begin{teo} The following statements are equivalent:

(1${_X}_t$) $F$ is $\mu$-constant.

(2${_X}_t$) $\partial G/\partial t \in \overline{\{\partial
G/\partial x_j\}}^{\dagger}$.

\vspace{0.2cm}

(3${_X}_t$) $v(J_X)=\{(t,x) \in \mathbb{C}\times\mathbb{C}^n \ | \
B_i(t,x)=0, \ i=1,\ldots,r\}=\mathbb{C}\times\{0\}$ near $(0,0)$,
where $B_1,\ldots,B_r$ are the maximal order minors of the Jacobian
matrix of $(f,\phi_t)$ (partial derivatives with respect to $x$,
only) and $J_X$
    is the ideal generated by them.
 \end{teo}

\noindent \textbf{Proof:}

(1$_{X_t}$) $\Leftrightarrow$ (3$_{X_t}$):

It is the proof of (1$_X$) $\Leftrightarrow$ (6$_X$) in Theorem
\ref{teo demonstracao}.

\vspace{0.2cm}

(1$_{X_t}$) $\Rightarrow$ (2$_{X_t}$):

Suppose $F$ is $\mu$-constant. By the Lê-Greuel formula, we have
that $\mu(f_t|_X)=\mu(X_t,0)+\mu(X_t\cap f_t^{-1}(0),0)$. Thus,
$\mu(X_t,0)$ and $\mu(X_t\cap f_t^{-1}(0),0)$ are constants. Then,
${A_f}_t$ holds for the pair $(\mathcal{X}_0,\mathbb{C}\times\{0\})$
(see \cite[Theorem 5.8]{Gaffney e Massey}). Therefore, $\partial
G/\partial t \in \overline{\{\partial G/\partial x_j\}}^{\dagger}$
(see \cite[Lemma 5.1]{Gaffney Specialization}).

\vspace{0.2cm}

(2$_{X_t}$) $\Rightarrow$ (1$_{X_t}$):

Suppose that $\partial G/\partial t \in \overline{\{\partial
G/\partial x_j\}}^{\dagger}$ then ${A_f}_t$ holds for the pair
$(\mathcal{X}_0,\mathbb{C}\times\{0\})$ (see \cite[Lemma
5.1]{Gaffney Specialization}). Hence, the Buchsbaum-Rim multiplicity
of the module $\{\partial G/\partial x_j\}$ is constant (see
\cite[Theorem 3.2]{Kleiman}). Thus, $F$ is $\mu$-constant (see
\cite[Lemma 3.3]{Kleiman}). \fim

\vspace{0.4cm}

\noindent\textbf{Acknowledgements}

\vspace{0.3cm}

We thank A. Miranda and T. Gaffney for the fruitful
discussion on the issue.

\bibliographystyle{amsplain}

\end{document}